\documentclass[journal]{IEEEtran}

\usepackage{graphicx}          % Include this line if your 
\usepackage{url}
\usepackage{amsmath,mathtools,amssymb}
\usepackage{bm}
\usepackage[inline]{enumitem}
\usepackage{color}
\usepackage[noadjust]{cite}
\usepackage{amsthm}

\newenvironment{smallarray}[1]
 {\null\,\vcenter\bgroup\scriptsize
  \arraycolsep=.13885em
  \hbox\bgroup$\array{@{}#1@{}}}
 {\endarray$\egroup\egroup\,\null}

\DeclareRobustCommand{\svdots}{% s for `scaling'
  \vbox{%
    \baselineskip=0.33333\normalbaselineskip
    \lineskiplimit=0pt
    \hbox{.}\hbox{.}\hbox{.}%
    \kern-0.2\baselineskip
  }%
}

\newtheorem{lem}{Lemma}
\newtheorem{thm}{Theorem}

\newtheorem{defn}{Definition}

\begin{document}

% \begin{frontmatter}

\title{Certificates of Nonexistence for Lyapunov-Based Stability, Stabilizability and Detectability of\\ LPV Systems} % Title, preferably not more 

% \thanks[footnoteinfo]{Co-funded by the European Union (ERC, Proacthis, no. 101055384). Views and opinions expressed are however those of the author(s) only and do not necessarily reflect those of the European Union or the European Research Council. Neither the European Union nor the granting authority can be held responsible for them.\\ Corresponding author: T.~J.~{Meijer}}

% \author[TUe]{T.~J.~Meijer}\ead{t.j.meijer@tue.nl},      % Add the 
% \author[ASML]{V.~S.~Dolk}\ead{victor.dolk@asml.com},    % e-mail address 
% \author[TUe]{W.~P.~M.~H.~Heemels}\ead{m.heemels@tue.nl} % (ead) as shown

% \address[TUe]{Department of Mechanical Engineering, Eindhoven University of Technology, Eindhoven, The Netherlands}  % Please supply                                              
% \address[ASML]{ASML, De Run 6665, 5504 DT Veldhoven, The Netherlands}             % full addresses here.
\author{T.~J.~{Meijer}, V.~S.~{Dolk} and W.~P.~M.~H.~{Heemels}
\thanks{Tomas Meijer and Maurice Heemels are with the Department of Mechanical Engineering, Eindhoven University of Technology, The Netherlands. (e-mail: t.j.meijer@tue.nl; m.heemels@tue.nl).}%
\thanks{Victor Dolk is with ASML, De Run 6665, 5504 DT Veldhoven, The Netherlands. (e-mail: victor.dolk@asml.com)}%
\thanks{Co-funded by the European Union (ERC, Proacthis, no. 101055384). Views and opinions expressed are however those of the author(s) only and do not necessarily reflect those of the European Union or the European Research Council. Neither the European Union nor the granting authority can be held responsible for them.}}

\maketitle
          
                             % keyword list or with the 
                                          % help of the Automatica 
                                          % keyword wizard

\begin{abstract}                          % Abstract of not more than 200 words.
	By computing Lyapunov functions of a certain, convenient structure, Lyapunov-based methods guarantee stability properties of the system or, when performing synthesis, of the relevant closed-loop or error dynamics. In doing so, they provide conclusive affirmative answers to many analysis and design questions in systems and control. When these methods fail to produce a feasible solution, however, they often remain inconclusive due to \begin{enumerate*}[label=(\alph*)] \item the method being conservative or \item the fact that there may be multiple causes for infeasibility, such as ill-conditioning, solver tolerances or true infeasibility.\end{enumerate*} To overcome this, we develop LMI-based theorems of alternatives based upon which we can guarantee, by computing a so-called certificate of nonexistence, that no poly-quadratic Lyapunov function exists for a given linear parameter-varying system. We extend these ideas to also certify the nonexistence of controllers and observers for which the corresponding closed-loop/error dynamics admit a poly-quadratic Lyapunov function. Finally, we illustrate our results in some numerical case studies.
\end{abstract}


% \end{frontmatter}

\section{Introduction}
Lyapunov's second method is arguably the most widely adopted technique for analyzing stability properties of dynamical systems~\cite{Kellett2015,Khalil1996-2nd-edition}. These methods work by computing a \textit{Lyapunov-function certificate}, i.e., by constructing a Lyapunov function (LF), that guarantees (asymptotic) stability of the system at hand. By applying a similar approach to the closed-loop or error dynamics for, respectively, some to-be-designed controller or observer, Lyapunov theory has also greatly impacted many synthesis problems. Over the years, many powerful computational approaches using, for instance, linear programs, see, e.g.,~\cite{Rantzer2018,Lazar2010}, or LMIs~\cite{Boyd1994,Scherer2000} have been developed to construct a wide variety of different classes of Lyapunov functions, such as, e.g., common quadratic LFs (cQLFs)~\cite{Boyd1994,Hespanha2018}, (piecewise) affine LFs~\cite{Milani2002,Rantzer2018}, piecewise quadratic LFs~\cite{Johansson1998,Decarlo2000}, poly-quadratic Lyapunov functions (poly-QLFs)~\cite{Daafouz2001,Heemels2010,Pandey2018,Meijer2023-lpv-arxiv} and sum-of-squares~\cite{Shivakumar2020,Parrilo2000}. While these methods are very powerful in applications where they produce a solution, they often remain \textit{inconclusive} when presented with a problem that is, for example, ill-posed or infeasible. Knowing whether the problem is truly infeasible is, however, also very valuable since it may indicate that we require a broader class of LFs or that we may have to redesign the underlying physical system. Nonetheless, very few computational approaches aimed towards providing these ``opposite'' results exist. The present paper aims to fill this gap.

In this paper, we consider a class of discrete-time (DT) polytopic linear parameter-varying (LPV) systems, see, e.g.,~\cite{Daafouz2001,Heemels2010}, which feature a specific mathematical structure that has enabled the development of many LMI-based analysis and synthesis conditions, see, e.g.,~\cite{Daafouz2001,Pandey2018,Oliveira1999,Meijer2023-lpv-arxiv} for stability analysis,~\cite{Meijer2023-lpv-arxiv,Daafouz2001,Pereira2021} for controller synthesis and~\cite{Meijer2023-lpv-arxiv,Heemels2010,Meijer2022-nl-obs-arxiv,Oliveira2022} for observer synthesis, that can be systematically solved. These methods are based on a class of parameter-dependent LFs known as poly-QLFs, which feature a quadratic state dependence and polytopic parameter dependence. The conditions introduced in~\cite{Daafouz2001,Heemels2010}--and, thereby, many conditions based on them--turn out to be not only sufficient but also necessary~\cite{Daafouz2001,Pandey2018}. Therefore, in theory, if we are unable to find a feasible solution to these conditions that should imply that no poly-QLF exists. In practice, however, there are many reasons that prevent the solver from finding a feasible solution, such as poor numerical conditioning, solver tolerances or true infeasibility, and, hence, often no conclusion can be drawn. To overcome this, we propose conditions that can be used to systematically construct so-called \textit{certificates of nonexistence} (CoNEs), i.e., numerical certificates that guarantee that a system does not admit a Lyapunov function of a certain structure. As we will see, these CoNEs are based on so-called theorems of alternatives (ToAs)~\cite{Balakrishnan2003}. To illustrate this, we use the following well-known ToA~\cite{Johansson1998,Liberzon2003} regarding cQLFs for continuous-time (CT) switched linear systems (SLS) (for the adopted notation, see also Section~\ref{sec:notation-bg} below).
\begin{thm}\label{thm:liberzon2003}
	For any matrices $A_i\in\mathbb{R}^{n\times n}$, $i\in\mathcal{N}\coloneqq\{1,2,\hdots,N\}$, $N\in\mathbb{N}_{\geqslant 1}$, exactly one of the following statements is true:
	\begin{enumerate}[labelindent=0pt,labelwidth=\widthof{\ref{item:liberzon2003-2}},label=(T\ref{thm:liberzon2003}.\arabic*),itemindent=0em,leftmargin=!]
		\item \label{item:liberzon2003-1} There exists a symmetric positive-definite matrix $P\in\mathbb{S}^{n}_{\succ 0}$ such that $A_i^\top P+PA_i \prec 0$ for all $i\in\mathcal{N}$.
		\item \label{item:liberzon2003-2} There exist symmetric positive-semidefinite matrices $R_0\in\mathbb{S}^{n}_{\succcurlyeq 0}$ and $R_i\in\mathbb{S}_{\succcurlyeq 0}^{n}$, $i\in\mathcal{N}$, such that $R_i\neq 0$, for some $i\in\{0\}\cup \mathcal{N}$, and 
		\begin{equation*}
			R_0 = \sum_{\mathclap{i\in\mathcal{N}}} A_iR_i+R_iA_i^\top.
		\end{equation*}
	\end{enumerate}
\end{thm}
\noindent Clearly,~\ref{item:liberzon2003-1} holds if and only if the CT SLS with matrices $\{A_i\}_{i\in\mathcal{N}}$, admits a cQLF~\cite{Liberzon2003}. By constructing a solution $\{R_i\}_{i\in\mathcal{N}}$ to~\ref{item:liberzon2003-2}, we certify that the system does not admit a cQLF, i.e., $\{R_i\}_{i\in\mathcal{N}}$ serves as a CoNE. 

We develop, in this paper, LMI-based conditions that can be used to compute CoNEs guaranteeing the nonexistence of a poly-QLF for a given DT polytopic system. We note that the use of poly-QLFs in the context of polytopic systems is equivalent to the use of mode-dependent Lyapunov functions for switched linear systems, see, e.g.,~\cite{Mason2012} and, thus, our results can also be applied to SLSs~\cite{Athanasopoulos2014,Philippe2016,Liberzon2003}. The development of our results is largely based on duality in semidefinite programming and theorems of alternatives, see, e.g.,~\cite{Balakrishnan2003}, which, in fact, can also be used to derive Theorem~\ref{thm:liberzon2003}. 

To illustrate the utility of constructing CoNEs, we consider (in Section~\ref{sec:example}) an example of a polytopic system that is GUAS but for which we show that no poly-QLF exists. Hence, the nonexistence of a poly-QLF for a polytopic system, in general, does not imply that the system is also not GUAS in contrast to what can be wrongly interpreted from~\cite[Theorem 1]{Daafouz2001}. We also develop LMI-based conditions to construct CoNEs for poly-quadratically (poly-Q) stabilizing controllers and observers that render the corresponding error system poly-Q stable (poly-QS), for which we also present a numerical example.

The remainder of this paper is organized as follows. We present some relevant notation and background in Section~\ref{sec:notation-bg}. Section~\ref{sec:cones} presents our main results, which we use in some applications and numerical examples in Section~\ref{sec:applications}. Finally, Section~\ref{sec:conclusions} presents some conclusions and the proofs of our results can be found in the Appendix.

\section{Notation and background}\label{sec:notation-bg}
\noindent{\bf Notation.} Let $\mathbb{R}=(-\infty,\infty)$ and $\mathbb{N}=\{0,1,2,\hdots\}$ denote the sets of real and non-negative natural numbers, respectively. We denote $\mathbb{R}_{\geqslant r}=[r,\infty)$ and $\mathbb{N}_{[n,m]}=\{n,n+1,\hdots,m\}$. The set of $n$-by-$n$ symmetric matrices is denoted by $\mathbb{S}^{n}$ and $M\prec 0$, $M\succ 0$, $M\preccurlyeq 0$ and $M\succcurlyeq 0$ mean, respectively, that $M\in\mathbb{S}^{n}$ is negative definite, positive definite, negative semidefinite and positive semidefinite. We denote, respectively, the set of such matrices by $\mathbb{S}^n_{\prec 0}$, $\mathbb{S}^{n}_{\succ 0}$, $\mathbb{S}^{n}_{\preccurlyeq 0}$ and $\mathbb{S}^{n}_{\succcurlyeq 0}$. The notation $(u,v)$ stand for $\left[\begin{smallarray}{cc}u^\top & v^\top\end{smallarray}\right]^\top$ and we use $\bm{e}_i\in\mathbb{R}^n$, $i\in\mathbb{N}_{[1,n]}$, to denote the $i$-th $n$-dimensional unit vector. For a set of matrices $\{A_i\}_{i\in\mathcal{N}}$, $\mathcal{N}\coloneqq\mathbb{N}_{[1,N]}$, $\operatorname{diag}\{A_1,A_2,\hdots,A_N\}$ denotes a block-diagonal matrix with $N$ diagonal blocks $A_i$, $i\in\mathcal{N}$. The joint spectral radius~\cite{Jungers2009} of a bounded set of matrices $\{A_i\}_{i\in\mathcal{N}}$ is denoted by $\rho(\{A_i\}_{i\in\mathcal{N}})\coloneqq\limsup_{t\rightarrow\infty}\{\|A_{i_1}\hdots A_{i_t}\|^{1/t}\mid i_t\in\mathcal{N}\text{ for all }t\in\mathbb{N}\}$. $\operatorname{Tr}(A)$ is the trace of $A\in\mathbb{R}^{n\times n}$, $n\in\mathbb{N}$. A continuous function $\alpha\colon\mathbb{R}_{\geqslant 0}\rightarrow\mathbb{R}_{\geqslant 0}$ is a $\mathcal{K}_{\infty}$-function, if it is strictly increasing, $\alpha\left(0\right)=0$ and $\alpha(r)\rightarrow\infty$ as $r\rightarrow\infty$. 

\subsection{Polytopic LPV systems}
We consider DT LPV systems of the form\footnote{We require $B$ and $C$ to be constant, for which we have necessary (and sufficient) LMI-based synthesis conditions, see, e.g.,~\cite{Meijer2023-lpv-arxiv}. This necessity plays an important role in developing our results. Recent works, see, e.g.,~\cite{Pereira2021,Oliveira2022}, present necessary and sufficient conditions, based on the so-called simplex approach~\cite{DeCaigny2010,Montagner2005}, for parameter-varying $B$ and $C$, which gives reason to believe that our results may be extended similarly.}
\begin{equation}
	\begin{aligned}
	x_{k+1} &= A(p_k)x_k + Bu_k,\\
	y_k &= Cx_k,
	\end{aligned}
	\label{eq:system}
\end{equation}
where $x_k\in\mathbb{R}^{n_x}$, $u_k\in\mathbb{R}^{n_u}$ and $y_k\in\mathbb{R}^{n_y}$ denote, respectively, the state, input and output at time $k\in\mathbb{N}$. The unknown parameter $p_k\in\mathbb{P}$, $k\in\mathbb{N}$, can be time-varying but remains within a known compact set $\mathbb{P}\subset\mathbb{R}^{n_p}$, as formalized below.
\begin{defn}
	The set of admissible parameter sequences $\mathcal{P}$ is the set of all parameter sequences $\bm{p}\coloneqq\{p_k\}_{k\in\mathbb{N}}$ with $p_k\in\mathbb{P}$ for all $k\in\mathbb{N}$.
\end{defn}
\noindent The matrix-valued function $A\colon\mathbb{P}\rightarrow\mathbb{R}^{n_x\times n_x}$ and matrices $B\in\mathbb{R}^{n_x\times n_u}$ and $C\in\mathbb{R}^{n_y\times n_x}$ are known. The class of \textit{polytopic} systems~\eqref{eq:system} is characterized by the following.
\begin{defn}\label{dfn:polytopic}
	The system~\eqref{eq:system} is said to be polytopic if there exist functions $\xi_i\colon\mathbb{P}\rightarrow\mathbb{R}_{\geqslant 0}$, $i\in\mathcal{N}\coloneqq\mathbb{N}_{[1,N]}$, such that the mapping $\xi\coloneqq(\xi_1,\xi_2,\hdots,\xi_N)$ satisfies $\xi(\mathbb{P})\subseteq \mathbb{X}\coloneqq\{\mu=(\mu_1,\mu_2,\hdots,\mu_N)\in\mathbb{R}^N_{\geqslant 0}\mid\sum_{i\in\mathcal{N}}\mu_i=1\}$, and matrices $A_i\in\mathbb{R}^{n_x\times n_x}$, $i\in\mathcal{N}$, such that, for all $\pi\in\mathbb{P}$,
	\begin{equation*}
		A(\pi) = \sum_{i\in\mathcal{N}}\xi_i(\pi)A_i.
	\end{equation*}
	If, in addition, the functions $\xi_i$, $i\in\mathcal{N}$, satisfy $\{\bm{e}_i\}_{i\in\mathcal{N}}\subset\xi(\mathbb{P})$, the system~\eqref{eq:system} is said to be strictly polytopic.
\end{defn}
\noindent Many works, such as, e.g., the aforementioned~\cite{Daafouz2001,Heemels2010,Pereira2021,Oliveira2022}, exploit the polytopic structure of the system~\eqref{eq:system} combined with a Lyapunov function with a similar polytopic structure, i.e., a poly-QLF as introduced in the next section, to formulate LMI-based analysis and synthesis conditions. Strict polytopicity essentially means that no conservatism is introduced to obtain a polytopic representation, i.e., for all vertices $\{\bm{e}_i\}_{i\in\mathcal{N}}$ of the (embedding) polytope, there exists $\nu_i\in\mathbb{P}$ such that $\bm{e}_i=\xi(\nu_i)$.

\subsection{Poly-quadratic stability}
Let $x_k(\bm{p},\bm{u},x_0)$ denote the state at time $k\in\mathbb{N}$ satisfying~\eqref{eq:system} with initial state $x_0\in\mathbb{R}^{n_x}$ at time $k=0$, parameter sequence $\bm{p}\in\mathcal{P}$ and input sequence $\bm{u}=\{u_j\}_{j\in\mathbb{N}}$ with $u_j\in\mathbb{R}^{n_u}$ for all $j\in\mathbb{N}$. For system~\eqref{eq:system} with $0$-input, i.e., $\bm{u}=\bm{0}\coloneqq \{0\}_{j\in\mathbb{N}}$, we define GUAS~\cite{Khalil1996-2nd-edition} below.
\begin{defn}\label{dfn:guas}
	The system~\eqref{eq:system} is said to be GUAS, if
	\begin{enumerate}[labelindent=0pt,labelwidth=\widthof{\ref{item:def-guas-2}},label=(D\ref{dfn:guas}.\arabic*),itemindent=0em,leftmargin=!]
		\item \label{item:def-guas-1} for any positive $\epsilon\in\mathbb{R}_{>0}$, there exists a positive $\delta(\epsilon)\in\mathbb{R}_{>0}$ (independent of $\bm{p}$) such that
		\begin{equation*}
			\|x_0\|<\delta(\epsilon) \implies \|x_k(\bm{p},\bm{0},x_0)\|<\epsilon,
		\end{equation*}
		for all $k\in\mathbb{N}_{\geqslant 1}$ and for all $\bm{p}\in\mathcal{P}$, and
		\item \label{item:def-guas-2} for any positive $\rho,\epsilon\in\mathbb{R}_{>0}$, there exists $K(\rho,\epsilon)\in\mathbb{N}$ (independent of $\bm{p}$) such that, for all $k>K(\rho,\epsilon)$ and for all $\bm{p}\in\mathcal{P}$,
		\begin{equation*}
			\|x_0\|<\rho \implies \|x_k(\bm{p},\bm{0},x_0)\|<\epsilon.
		\end{equation*}
	\end{enumerate}
\end{defn}
\noindent The definition above, in which we consider uniformity with respect to $\bm{p}$ (and $x_0$), is obtained by applying the general definition of GUAS, see, e.g.,~\cite{Khalil1996-2nd-edition}, to the class of LPV systems considered here. Since it is generally difficult to verify whether a given (polytopic) system~\eqref{eq:system} is GUAS in the sense of Definition~\ref{dfn:guas}, much research has focussed on proposing Lyapunov-based techniques that, often at the cost of introducing some degree of conservatism, can be systematically verified. As mentioned before, the specific class of parameter-dependent quadratic Lyapunov functions known as poly-QLFs, introduced below, is often used for polytopic systems~\eqref{eq:system} because they \begin{enumerate*}[label=(\alph*)] \item can be constructed using LMI-based techniques, and \item are less conservative than traditional cQLFs\end{enumerate*}~\cite{Daafouz2001}. We establish the corresponding notion of poly-QS below.
\begin{defn}
    The system~\eqref{eq:system}, being polytopic, is said to be poly-QS, if it admits a poly-QLF, i.e., a function $V\colon\mathbb{P}\times\mathbb{R}^{n_x}\rightarrow\mathbb{R}_{\geqslant 0}$ of the form
	\begin{equation}
		V(\pi,x) = x^\top \sum_{i\in\mathcal{N}}\xi_i(\pi)P_ix,
		\label{eq:poly-lyap}
	\end{equation}
    where $P_i\in\mathbb{S}^{n_x}$, $i\in\mathcal{N}$, and $\xi_i$, $i\in\mathcal{N}$, are the functions in Definition~\ref{dfn:polytopic}, satisfying, for some $a_j\in\mathbb{R}_{>0}$, $j\in\{1,2,3\}$,
	\begin{align}
		a_1\|x\|^2 &\leqslant V(\pi,x)\leqslant a_2\|x\|^2,\label{eq:sw-bnds}\\
		V(\pi_+,x_+) &\leqslant V(\pi,x)-a_3\|x\|^2,\label{eq:desc-cond}
	\end{align} 
    for all $\pi_+,\pi\in\mathbb{P}$ and $x_+,x\in\mathbb{R}^{n_x}$ with $x_+=A(\pi)x$.
\end{defn}
\noindent Poly-QS implies GUAS~\cite{Khalil1996-2nd-edition} (and, in fact, even uniform global exponential stability (UGES)~\cite{Jiang2002}), however, the converse is not true, as we will also illustrate by means of an example in Section~\ref{sec:applications}. In the remainder, we will develop conditions that, if they are true for a given system~\eqref{eq:system}, provide CoNEs guaranteeing that the system does not admit a poly-QLF and, thereby, providing the possibility to conclusively show that a system is not poly-QS. We further build on these CoNEs to develop conditions that guarantee the nonexistence of poly-Q stabilizing controllers and observers that render the corresponding error system poly-QS. Thereby, we show the application of CoNEs to synthesis problems as well (and thus certain Lyapunov-based stabilizability and detectability conditions).

\section{Certificates of nonexistence}\label{sec:cones}
\subsection{Stability analysis}
To facilitate the computation of CoNEs for poly-QLFs, we provide a ToA below, of which a proof can be found in the Appendix.
\begin{thm}\label{thm:toa-poly-QLF}
	Consider the system~\eqref{eq:system}, being strictly polytopic. Exactly one of the following statements is true:
	\begin{enumerate}[labelindent=0pt,labelwidth=\widthof{\ref{item:toa-poly-QLF-2}},label=(T\ref{thm:toa-poly-QLF}.\arabic*),itemindent=0em,leftmargin=!]
		\item \label{item:toa-poly-QLF-1} The system~\eqref{eq:system} is poly-QS.
		\item \label{item:toa-poly-QLF-2} There exist symmetric positive-semidefinite matrices $R_i\in\mathbb{S}^{n_x}_{\succcurlyeq 0}$, $i\in\mathbb{N}_{[1,N^2]}$, such that
		\begin{equation*}
			\sum_{j\in\mathcal{N}} A_jR_{N(j-1)+i}A_j^\top - R_{N(i-1)+j}\succcurlyeq 0,
		\end{equation*}
		for all $i\in\mathcal{N}$, and $R_i\neq 0$ for some $i\in\mathbb{N}_{[1,N^2]}$.
	\end{enumerate}
\end{thm}
\noindent It follows from Theorem~\ref{thm:toa-poly-QLF} that, by finding $R_i$, $i\in\mathbb{N}_{[1,N^2]}$, satisfying the conditions in~\ref{item:toa-poly-QLF-2}, we can conclude that the system at hand does not admit a poly-QLF. This is a powerful result since it allows us to give a conclusive ``negative'' answer to the question ``Is the strictly polytopic system~\eqref{eq:system} poly-QS?'' whereas existing conditions to construct a poly-QLF are often only capable of providing conclusive affirmative (but not ``negative'') answers. Moreover, a solution $R_i$, $i\in\mathbb{N}_{[1,N^2]}$, to~\ref{item:toa-poly-QLF-2}, can then be presented as concrete numerical evidence for such a ``negative'' result.

\subsection{Observer design}
Next, we extend consider CoNEs in the context of designing linear parameter-varying observers of the form
\begin{equation*}
	\hat{x}_{k+1} = A(p_k)\hat{x}_k + Bu_k + L(k,\bm{p})(C\hat{x}_k-y_k),
\end{equation*}
where $\hat{x}_k\in\mathbb{R}^{n_x}$ denotes the estimated state at time $k\in\mathbb{N}$. We aim to design the observer gain $L\colon\mathbb{N}\times\mathcal{P}\rightarrow\mathbb{R}^{n_x\times n_x}$, which depends on the entire parameter sequence $\bm{p}$, such that the estimation error system
\begin{equation}
	e_{k+1} = (A(p_k) + L(k,\bm{p})C)e_k,
	\label{eq:error-system}
\end{equation}
with $e_k\coloneqq \hat{x}_k-x_k$, is GUAS. In particular, we focus on gains $L$ that render~\eqref{eq:error-system} poly-QS, as formalized below.
\begin{defn}
	The system~\eqref{eq:system}, being polytopic, is said to be poly-Q detectable, if there exists a matrix-valued function $L\colon\mathbb{N}\times\mathcal{P}\rightarrow\mathbb{R}^{n_x\times n_y}$ that renders the estimation error system~\eqref{eq:error-system} poly-QS.
\end{defn}

To check whether a system is poly-Q detectable or not, we can use the result below, which provides a necessary and sufficient condition and, hence, depending on whether it admits a solution, we should in theory be able to use it to either affirm or disprove poly-Q detectability.
\begin{thm}{\cite[Theorem 1]{Meijer2023-lpv-arxiv}}\label{thm:pQdet}
	The system~\eqref{eq:system}, being strictly polytopic, is poly-Q detectable if and only if there exist symmetric positive-definite matrices $P_i\in\mathbb{S}^{n_x}_{\succ 0}$, $i\in\mathcal{N}$, such that, for all $i,j\in\mathcal{N}$,
	\begin{equation}
		P_i-A_i^\top P_jA_i +C^\top C\succ 0.
		\label{eq:polyQdetcond}
		\vspace*{-\belowdisplayskip}
	\end{equation}
\end{thm}
\noindent If a solution to the conditions in Theorem~\ref{thm:pQdet} is found or, alternatively, if an observer gain $L$ that renders~\eqref{eq:error-system} poly-QS is found using, for instance, the synthesis conditions in~\cite{Heemels2010,Pandey2018}, we can decidedly conclude that the system is poly-Q detectable. When we cannot find a feasible solution, the conditions in Theorem~\ref{thm:pQdet} we should, in theory, be able to conclude that the system is not poly-Q detectable, however, in practice, this is only one of the possible explanations along with, e.g., numerical conditioning or solver tolerances. The ToA below provides a set of alternative conditions to disprove poly-Q detectability.
\begin{thm}\label{thm:toa-pqdet}
	Consider the system~\eqref{eq:system}, being strictly polytopic. Exactly one of the following statements is true:
	\begin{enumerate}[labelindent=0pt,labelwidth=\widthof{\ref{item:toa-pqdet-2}},label=(T\ref{thm:toa-pqdet}.\arabic*),itemindent=0em,leftmargin=!]
		\item \label{item:toa-pqdet-1} The system~\eqref{eq:system} is poly-Q detectable.
		\item \label{item:toa-pqdet-2} There exist symmetric positive-semidefinite matrices $R_i\in\mathbb{S}^{n_x}_{\succcurlyeq 0}$, $i\in\mathbb{N}_{[1,N^2]}$, such that
		\begin{equation*}
				\sum_{j\in\mathcal{N}} A_jR_{N(j-1)+i}A_j^\top -R_{N(i-1)+j} \succcurlyeq 0,
		\end{equation*}
		for all $i\in\mathcal{N}$, $R_i\neq 0$ for some $i\in\mathbb{N}_{[1,N^2]}$ and
		\begin{equation*}
			C\sum_{\mathclap{i\in\mathbb{N}_{[1,N^2]}}}R_iC^\top = 0.
		\end{equation*}
	\end{enumerate}
\end{thm}
\noindent As before, any set of matrices $\{R_i\}_{i\in\mathbb{N}_{[1,N^2]}}$, satisfying~\ref{item:toa-pqdet-2}, serves as a CoNE for an observer gain $L$ that renders~\eqref{eq:error-system} poly-QS. Hence, by solving the LMI-based conditions in~\ref{item:toa-pqdet-2}, we are able to provide conclusive ``negative'' answers to the question ``Is the strictly polytopic system~\eqref{eq:system} poly-Q detectable?''.

\subsection{Controller design}
In this section, we consider CoNEs in the context of state-feedback controller design. To be precise, we consider the design of an LPV control law
\begin{equation}
    u_k = K(k,\bm{p})x_k,
    \label{eq:c-law}
\end{equation}
where $K\colon\mathbb{N}\times\mathcal{P}\rightarrow\mathbb{R}^{n_u\times n_x}$ is the to-be-designed controller gain, which is allowed to depend on the entire parameter sequence $\bm{p}$. We aim to design $K$ such that it renders the corresponding closed-loop system, given by
\begin{equation}
    x_{k+1} = (A(p_k)+BK(k,\bm{p}))x_k,
    \label{eq:closed-loop}
\end{equation}
poly-QS, as formalized below.
\begin{defn}\label{dfn:poly-q-stabilizability}
    The system~\eqref{eq:system}, being polytopic, is said to be poly-Q stabilizable, if there exists a matrix-valued function $K\colon\mathbb{N}\times\mathcal{P}\rightarrow\mathbb{R}^{n_u\times n_x}$ that renders the closed-loop system~\eqref{eq:closed-loop} poly-QS.
\end{defn}
Poly-Q stabilizability implies that some control law~\eqref{eq:c-law} renders~\eqref{eq:closed-loop} GUAS. Unlike for poly-Q detectability, to the best of the authors' knowledge, there does not exist a necessary and sufficient condition for poly-Q stabilizability. There are necessary and sufficient conditions for the existence of a poly-Q stabilizing $K$ with a predetermined polytopic structure, see, e.g.,~\cite{Daafouz2001}, however, not for the general $K$ as in Definition~\ref{dfn:poly-q-stabilizability}. There is, however, the necessary condition below, which is not sufficient (even for $K$ with a predetermined polytopic structure). 
\begin{thm}[{\cite[Theorem 2]{Meijer2023-lpv-arxiv}}]\label{thm:pQstab}
    Suppose the system~\eqref{eq:system}, being strictly polytopic, is poly-Q stabilizable. Then, there exist symmetric postive-definite matrices $S_i\in\mathbb{S}^{n_x}_{\succ 0}$, $i\in\mathcal{N}$, such that, for all $i,j\in\mathcal{N}$,
    \begin{equation}
        S_j - A_iS_iA_i^\top + BB^\top\succ 0.
        \label{eq:pqstab-cond}
    \end{equation}
\end{thm}
In contrast with Theorem~\ref{thm:pQdet}, we cannot conclude that the system~\eqref{eq:system} is poly-Q stabilizable by finding a solution $\{S_i\}_{i\in\mathcal{N}}$ to the conditions in Theorem~\ref{thm:pQstab}. Nevertheless, we can formulate the ToA below for the conditions in~\eqref{eq:pqstab-cond}.
\begin{lem}\label{lem:interm-toa}
    Consider the system~\eqref{eq:system}. Exactly one of the following statements is true:
    \begin{enumerate}[labelindent=0pt,labelwidth=\widthof{\ref{item:interm-toa-2}},label=(L\ref{lem:interm-toa}.\arabic*),itemindent=0em,leftmargin=!]
		\item \label{item:interm-toa-1} There exist symmetric positive-definite matrices $S_i\in\mathbb{S}^{n_x}_{\succ 0}$, $i\in\mathcal{N}$, satisfying~\eqref{eq:pqstab-cond} for all $i,j\in\mathcal{N}$.
		\item \label{item:interm-toa-2} There exist symmetric positive-semidefinite matrices $R_i\in\mathbb{S}^{n_x}_{\succcurlyeq 0}$, $i\in\mathbb{N}_{[1,N^2]}$, such that
		\begin{equation}\label{eq:interm-toa-2-a}
				\sum_{j\in\mathcal{N}} A_i^\top R_{N(j-1)+i}A_i -R_{N(i-1)+j} \succcurlyeq 0,
		\end{equation}
		for all $i\in\mathcal{N}$, $R_i\neq 0$ for some $i\in\mathbb{N}_{[1,N^2]}$ and
		\begin{equation}\label{eq:interm-toa-2-b}
			B^\top\sum_{\mathclap{i\in\mathbb{N}_{[1,N^2]}}}R_iB = 0.
		\end{equation}
	\end{enumerate}
\end{lem}
Since the conditions in Theorem~\ref{thm:pQstab} do not imply poly-Q stabilizability, we cannot conclude that the system is poly-Q stabilizable if~\ref{item:interm-toa-1} holds. However, we can still conclude that it is not poly-Q stabilizable if~\ref{item:interm-toa-2} holds.
\begin{thm}\label{thm:toa-pqstab}
    Consider the system~\eqref{eq:system}, being strictly polytopic. At most one of the following statements is true:
	\begin{enumerate}[labelindent=0pt,labelwidth=\widthof{\ref{item:toa-pqstab-2}},label=(T\ref{thm:toa-pqstab}.\arabic*),itemindent=0em,leftmargin=!]
		\item \label{item:toa-pqstab-1} The system~\eqref{eq:system} is poly-Q stabilizable.
		\item \label{item:toa-pqstab-2} There exist symmetric positive-semidefinite matrices $R_i\in\mathbb{S}^{n_x}_{\succcurlyeq 0}$, $i\in\mathbb{N}_{[1,N^2]}$, such that~\eqref{eq:interm-toa-2-a}-\eqref{eq:interm-toa-2-b} hold and $R_i\neq 0$, for some $i\in\mathbb{N}_{[1,N^2]}$.
	\end{enumerate}
\end{thm}
Unlike the ToAs presented before for poly-QS (Theorem~\ref{thm:toa-poly-QLF}) and poly-Q detectability (Theorem~\ref{thm:toa-pqdet}) as well as Lemma~\ref{lem:interm-toa}, which were strong ToAs in the sense that exactly one of their statements was true, Theorem~\ref{thm:toa-pqstab} is a so-called weak ToA, i.e., at most one of the statements is true. Regardless, any solution $\{R_i\}_{i\in\mathbb{N}_{[1,N^2]}}$ to~\ref{item:toa-pqstab-2} serves as a CoNE, which guarantees that no poly-Q stabilizing controller~\eqref{eq:c-law} exists. Theorem~\ref{thm:toa-pqstab} may, however, be conservative in the sense that there may exist systems that are not poly-Q stabilizable, but for which~\ref{item:toa-pqstab-2} is not feasible. In such cases, Theorem~\ref{thm:toa-pqstab} cannot be used to compute a CoNE. 

\section{Applications}\label{sec:applications}
\subsection{The inequivalence of poly-QS and GUAS}\label{sec:example}
It can be wrongly interpreted from~\cite[Theorem 1]{Daafouz2001} that the system~\eqref{eq:system}, being strictly polytopic, is GUAS if and only if it admits a poly-QLF.\footnote{Strictly polytopic systems~\eqref{eq:system} admit an LF of the form in~\cite[Theorem 1]{Daafouz2001} if and only if they admit a poly-QLF.} The sufficiency, indeed, follows from the fact that $V$ is a LF for the system~\eqref{eq:system}~\cite{Khalil1996-2nd-edition}. However, necessity does not hold as we will demonstrate next by means of an example, i.e., a strictly polytopic system~\eqref{eq:system} that is GUAS but does not admit a poly-QLF.

To this end, consider the system~\eqref{eq:system} with\footnote{The authors would like to express their gratitude to Rapha\"{e}l Jungers, Virginie Debauche and Matteo Della Rossa for their assistance in constructing this example.}
\begin{equation}
	A(\pi) = (1-\pi)A_1+\pi A_2,
	\label{eq:cex-A}
\end{equation}
where $p_k\in\mathbb{P}\coloneqq\left[0,1\right]$, which renders the system strictly polytopic, and
\begin{equation}
	A_1 = \epsilon\begin{bmatrix}
		0.80 & 0.65\\
		-0.34 & 0.90
	\end{bmatrix},\quad A_2 = \epsilon\begin{bmatrix}
		0.43 & 0.62\\
		-1.48 & 0.14
	\end{bmatrix},
	\label{eq:cex-mat}
\end{equation}
with $\epsilon = 1/1.19$. The system is GUAS if and only if the joint spectral radius (JSR) is strictly negative, i.e., $\rho(\{A_0,A_1\})<1$,~\cite[Corollary 1.1]{Jungers2009}. Using the JSR toolbox for MATLAB~\cite{Vankeerberghen2014}, we compute that $\rho(\{A_1,A_2\})< 0.9786$ and, hence, the system is GUAS. However, using the YALMIP toolbox~\cite{Lofberg2004} with MOSEK~\cite{Mosek2019}, we find that the matrices
\begin{equation*}
	\begin{aligned}
		R_1&=\begin{bmatrix}
			0.7 & 0.3\\
			0.3 & 2.3
		\end{bmatrix},\,R_2 = \begin{bmatrix}
			1.0 & 0.5\\
			0.5 & 1.5
		\end{bmatrix},\\
		R_3&= \begin{bmatrix}
			1.1 & -0.1\\
			-0.1 & 1.0
		\end{bmatrix},\text{ and }R_4=\begin{bmatrix}
			0.9 & 0.2\\
			0.2 & 1.2
		\end{bmatrix},
	\end{aligned}
\end{equation*}
satisfy~\ref{item:toa-poly-QLF-2} and, thus, constitute a CoNE for a poly-QLF for the considered system. Thus, the system~\eqref{eq:system} with~\eqref{eq:cex-A}-\eqref{eq:cex-mat} is GUAS, but not poly-QS.

\subsection{Poly-quadratic stabilizability}
Next, we consider poly-quadratic stabilizability of the strictly polytopic system~\eqref{eq:system} with $p_k\in\mathbb{P}\coloneqq [0,1]$, $A(\pi)=(1-\pi)\bar{A}_1+\pi \bar{A}_2$, $B=(0,1,0,0)$,
\begin{align*}
    \bar{A}_1 &= \begin{bmatrix}\begin{smallmatrix}
        A_1 & I\\
        0 & A_1
    \end{smallmatrix}\end{bmatrix}\text{ and }
    A_2=\begin{bmatrix}\begin{smallmatrix}
        A_2 & I\\
        0 & A_2
    \end{smallmatrix}\end{bmatrix},
\end{align*}
where $A_1$ and $A_2$ are as in~\eqref{eq:cex-mat}. For this system, 
\begin{equation*}
    \begin{aligned}
        &R_1 = \begin{bmatrix}\begin{smallmatrix}
            0    &        0     &       0  &          0\\
            0     &       0    &        0   &         0\\
            0      &      0   &1.65 &  2.23\\
            0       &     0  & 2.23  & 3.04
        \end{smallmatrix}\end{bmatrix},R_2 = \begin{bmatrix}\begin{smallmatrix}
            0     &       0        &    0    &        0\\
            0      &      0     &       0   &         0\\
            0      &      0 &  0.01 &  -0.02\\
            0   &         0 & -0.02 &   2.47
        \end{smallmatrix}\end{bmatrix},\\
        &R_3 = \begin{bmatrix}\begin{smallmatrix}
            0      &      0    &        0      &      0\\
            0    &        0    &        0     &       0\\
            0       &    0&   3.69  &0.38\\
            0     &       0 &  0.38  & 0.05
        \end{smallmatrix}\end{bmatrix},R_4 = \begin{bmatrix}\begin{smallmatrix}
            0      &      0       &     0        &    0\\
            0      &      0       &     0       &     0\\
            0       &     0   & 0.85 &  -0.37\\
            0           & 0 & -0.37 &   0.17

        \end{smallmatrix}\end{bmatrix},
    \end{aligned}
\end{equation*}
satisfy the conditions in~\ref{item:toa-pqstab-2} and, thereby, constitute a CoNE guaranteeing that no controller~\eqref{eq:c-law} exists that renders the corresponding closed-loop system~\eqref{eq:closed-loop} poly-QS. This does not mean that no controller exists that renders the closed-loop system~\eqref{eq:closed-loop} GUAS, just that no controller of the form~\eqref{eq:c-law} renders~\eqref{eq:closed-loop} poly-QS. In fact, we have $\rho(\{A_1,A_2\})<0.9786$ and, thus, the controller~\eqref{eq:c-law} with $K(k,\bm{p})=0$, $k\in\mathbb{N}$, $\bm{p}\in\mathcal{P}$, renders the corresponding closed-loop system~\eqref{eq:closed-loop} is GUAS.

\section{Conclusions}\label{sec:conclusions}
We presented, for a class of polytopic LPV systems, theorems of alternatives regarding the existence of poly-QLFs as well as poly-Q stabilizing controllers and observers, that yield poly-QS error dynamics. Moreover, we showed how these ToAs can be used to compute CoNEs, which can guarantee that such LFs, controllers and observers do not exist for a given system. Hence, in contrast to existing Lyapunov-based techniques, these CoNEs provide a powerful way of providing conclusive and verifiable ``negative'' answers to analysis and synthesis questions for polytopic systems. We illustrated the application of these CoNEs in providing an example, which establishes the inequivalence of poly-QS and GUAS, and analyzing poly-Q stabilizability. Future work should be directed towards extending these results to different classes of systems and LFs. Finally, it is interesting to investigate whether CoNEs can be used to derive insight into whether a specific, broader class of LFs can guarantee stability or how we can redesign our system to make it poly-Q detectable or stabilizable.

\bibliographystyle{plain}        % Include this if you use bibtex 
\bibliography{phd-bibtex}           % and a bib file to produce the 
                                 % bibliography (preferred). The
                                 % correct style is generated by
                                 % Elsevier at the time of printing.

\appendix
\section{Preliminaries}    % Each appendix must have a short title.
For the sake of self-containedness, we state the ToA below, which we use in proving our results.
\begin{thm}[{\cite[Theorem 1]{Balakrishnan2003}}]\label{thm:thm1balakrishnan}
	Let $\mathcal{S}$ be a space of block-diagonal Hermitian matrices with inner product 
	\begin{equation}\label{eq:inner-prod}
		\left\langle\operatorname{diag}\{A_i\}_{i\in\mathcal{N}},\operatorname{diag}\{B_j\}_{j\in\mathcal{N}}\right\rangle_{\mathcal{S}}=\sum_{k\in\mathcal{N}} \operatorname{Tr}A_kB_k,
	\end{equation}
	with $\mathcal{N}\coloneqq\mathbb{N}_{[1,N]}$, and suppose $\mathcal{V}$ is a finite-dimensional vector space with an inner product $\left\langle\cdot,\cdot\right\rangle_{\mathcal{V}}$, $\mathcal{A}\colon\mathcal{V}\rightarrow\mathcal{S}$ is a linear mapping and $A_0\in\mathcal{S}$. Moreover, let $\mathcal{A}^{\operatorname{adj}}\colon\mathcal{S}\rightarrow\mathcal{V}$ denote the adjoint mapping of $\mathcal{A}$, i.e., $\left\langle\mathcal{A}(x),R\right\rangle_{\mathcal{S}}=\left\langle x,\mathcal{A}^{\operatorname{adj}}(R)\right\rangle_{\mathcal{V}}$ for all $x\in\mathcal{V}$ and $R\in\mathcal{S}$. Then, exactly one of the following statements is true:
	\begin{enumerate}[labelindent=0pt,labelwidth=\widthof{\ref{item:balakrishnan-2}},label=(T\ref{thm:thm1balakrishnan}.\arabic*),itemindent=0em,leftmargin=!]
		\item \label{item:balakrishnan-1} There exists an $x\in\mathcal{V}$ with $\mathcal{A}(x)+A_0\succ 0$.
		\item \label{item:balakrishnan-2} There exists an $R\in\mathcal{S}$ with $R\succcurlyeq 0$, $R\neq 0$, $\mathcal{A}^{\operatorname{adj}}(R)=0$ and $\left\langle A_0,R\right\rangle_{\mathcal{S}}\leqslant 0$.
	\end{enumerate}
\end{thm}
In Theorem~\ref{thm:thm1balakrishnan}, $A\succ 0$ and $A\succcurlyeq 0$ still denote positive definiteness and semidefiniteness in the sense that the smallest eigenvalue of $A$ is, respectively, positive and non-negative (not with respect to the inner product~\eqref{eq:inner-prod}).

\section{Proofs}         % Sections and subsections are supported  
                                        % in the appendices.
\begin{proof}[Proof of Theorem~\ref{thm:toa-poly-QLF}]
	Consider the system~\eqref{eq:system}, being strictly polytopic. Using its strict polytopicity, the system~\eqref{eq:system} is poly-QS if and only if there exist symmetric matrices $P_i\in\mathbb{S}^{n_x}$, $i\in\mathcal{N}$, satisfying, for all $i,j\in\mathcal{N}$,~\cite{Daafouz2001}
	\begin{equation}
		P_i\succ 0,\text{ and }P_i - A_i^\top P_jA_i\succ 0,\text{ for all }i,j\in\mathcal{N}.
		\label{eq:poly-QScond}
	\end{equation}	
	The remainder of our proof of Theorem~\ref{thm:toa-poly-QLF} is based on Theorem~\ref{thm:thm1balakrishnan}. Let $\mathcal{S}^{n}_N\subseteq\mathbb{S}^{Nn}$ denote a space of block-diagonal Hermitian matrices, composed of $N$ diagonal blocks of size $n$-by-$n$, with inner product~\eqref{eq:inner-prod} and let
	\begin{equation}
		P\coloneqq\operatorname{diag}\{P_1,P_2,\hdots,P_N\}.
		\label{eq:Pdiag}
	\end{equation}
	Then,~\eqref{eq:poly-QScond} can be expressed as $\mathcal{A}(P)\succ 0$ with $\mathcal{A}\colon\mathcal{S}^{n_x}_N\rightarrow\mathcal{S}^{n_x}_{(N+1)N}$ (i.e., $\mathcal{V}\leftarrow\mathcal{S}^{n_x}_N$ and $\mathcal{S}\leftarrow\mathcal{S}^{n_x}_{(N+1)N}$ in Theorem~\ref{thm:thm1balakrishnan}) given by
	\begin{equation}
		\mathcal{A}(P)\coloneqq \operatorname{diag}\{\operatorname{diag}\{P_i-A_i^\top P_jA_i\}_{i,j\in\mathcal{N}},P\}.
		\label{eq:Acal}
	\end{equation}
    Note that we use the convention $\{A_{ij}\}_{i,j\in\mathcal{N}}=\{A_{11},A_{21},A_{12},A_{22}\}$. Using the fact that $\operatorname{Tr}ABC=\operatorname{Tr}BCA$ and that $\mathcal{A}^{\operatorname{adj}}\colon\mathcal{S}^{n_x}_{(N+1)N}\rightarrow\mathcal{S}^{n_x}_{N}$ satisfies, for all $P_i\in\mathbb{S}^{n_x}$, $i\in\mathcal{N}$, and $R_i\in\mathbb{S}^{n_x}$, $i\in\mathbb{N}_{[1,(N+1)N]}$, 
	\begin{align}
			&\left\langle\mathcal{A}(P),R\right\rangle_{\mathcal{S}^{n_x}_{(N+1)N}}=\sum_{\mathclap{i\in\mathcal{N}}} \operatorname{Tr}P_iR_{N^2+i} + \nonumber\\
			&\, \sum_{\mathclap{i\in\mathcal{N}}}\sum_{\mathclap{j\in\mathcal{N}}} \operatorname{Tr}(P_i-A_i^\top P_jA_i)R_{(j-1)N+i}=\nonumber\\
			& \sum_{i\in\mathcal{N}}\operatorname{Tr}P_i(R_{N^2+i}+\sum_{j\in\mathcal{N}}(R_{(j-1)N+i}-A_jR_{(i-1)N+j}A_j^\top))\nonumber\\
			&=\left\langle P,\mathcal{A}^{\operatorname{adj}}(R)\right\rangle_{\mathcal{S}^{n_x}_{N}}\label{eq:Aadj-deriv}
	\end{align}
	with $R=\operatorname{diag}\{R_i\}_{i\in\mathbb{N}_{[1,(N+1)N]}}\in\mathcal{S}^{n_x}_{(N+1)N}$, we find
	\begin{align}
		&\mathcal{A}^{\operatorname{adj}}(R) = \label{eq:Aadj}\\
		&\operatorname{diag}\{R_{N^2+i}+\sum_{\mathclap{j\in\mathcal{N}}} (R_{(j-1)N+i}-A_jR_{(i-1)N+j}A_j^\top)\}_{i\in\mathcal{N}}.\nonumber
	\end{align}
	By Theorem~\ref{thm:thm1balakrishnan} (with $A_0=0$), there exist $P_i\in\mathbb{S}^{n_x}$, $i\in\mathcal{N}$, satisfying~\eqref{eq:poly-QScond} for all $i,j\in\mathcal{N}$, if and only if there does not exist $R=\operatorname{diag}\{R_i\}_{i\in\mathbb{N}_{[1,(N+1)N]}}\in\mathcal{S}^{n_x}_{(N+1)N}$ such that $R\succcurlyeq 0$, $R\neq 0$ and $\mathcal{A}^{\operatorname{adj}}(R)=0$, where $\mathcal{A}^{\operatorname{adj}}(R)=0$ evaluates to, for all $i\in\mathcal{N}$,
	\begin{equation}
		R_{N^2+i}=\sum_{\mathclap{j\in\mathcal{N}}}A_jR_{(i-1)N+j}A_j^\top - R_{(j-1)N+i}.
		\label{eq:lastRs}
	\end{equation}
	Using that $R_{N^2+i}\succcurlyeq 0$ for all $i\in\mathcal{N}$ and that at least one of the first $N^2$ $R_i$'s is nonzero, i.e., $R_i\neq 0$ for some $i\in\mathbb{N}_{[1,N^2]}$ (otherwise $R=0$ by~\eqref{eq:lastRs}), we conclude that exactly one of the statements in Theorem~\ref{thm:toa-poly-QLF} holds.
\end{proof}

\begin{proof}{Proof of Theorem~\ref{thm:toa-pqdet}}
	Consider the system~\eqref{eq:system}, being strictly polytopic. By Theorem~\ref{thm:toa-pqdet}, the system~\eqref{eq:system} is poly-Q detectable if and only if there exist symmetric positive-definite matrices $P_i\in\mathbb{S}^{n_x}_{\succ 0}$, $i\in\mathcal{N}$, such that~\eqref{eq:polyQdetcond} holds for all $i,j\in\mathcal{N}$.	The remainder of our proof of Theorem~\ref{thm:toa-pqdet} is based on Theorem~\ref{thm:thm1balakrishnan}. Let $\mathcal{S}^{n}_N\subseteq\mathbb{S}^{Nn}$ denote, again, a space of block-diagonal Hermitian matrices, composed of $N$ diagonal blocks of size $n$-by-$n$, with inner product~\eqref{eq:inner-prod} and let $P$ be as in~\eqref{eq:Pdiag} such that we can express~\eqref{eq:polyQdetcond} as $\mathcal{A}(P)+A_0\succ 0$ with~\eqref{eq:Acal} and $A_0\coloneqq\operatorname{diag}\{\operatorname{diag}\{C^\top C\}_{i,j\in\mathcal{N}},0_{Nn_x\times Nn_x}\}$. Since $\mathcal{A}$ is the same as in~\eqref{eq:Acal}, we can use~\eqref{eq:Aadj-deriv} to find that $A^{\operatorname{adj}}\colon\mathcal{S}^{n_x}_{(N+1)N}\rightarrow\mathcal{S}^{n_x}_{N}$ is as in~\eqref{eq:Aadj}. Hence, the main difference with respect to the proof of Theorem~\ref{thm:toa-poly-QLF} is the presence of a nonzero $A_0$. 
	
	We apply Theorem~\ref{thm:thm1balakrishnan} to find that there exist $P_i\in\mathbb{S}^{n_x}$, $i\in\mathcal{N}$, satisfying~\eqref{eq:polyQdetcond} if and only if there does not exist $R=\operatorname{diag}\{R_i\}_{i\in\mathbb{N}_{[1,(N+1)N]}}\in\mathcal{S}^{n_x}_{(N+1)N}$ such that $R\succcurlyeq 0$, $R\neq 0$, $\mathcal{A}^{\operatorname{adj}}(R)=0$ and $\left\langle A_0,R\right\rangle_{\mathcal{S}^{n_x}_{(N+1)N}}\leqslant 0$, where $\mathcal{A}^{\operatorname{adj}}(R)=0$ evaluates to~\eqref{eq:lastRs}, for all $i\in\mathcal{N}$, and 
	\begin{equation}
		\left\langle A_0,R\right\rangle_{\mathcal{S}^{n_x}_{(N+1)N}} = \operatorname{Tr}C\sum_{\mathclap{i\in\mathbb{N}_{[1,N^2]}}}R_iC^\top\leqslant 0,
		\label{eq:negsemdeftr}
	\end{equation}
	where we used that $\operatorname{Tr}AB = \operatorname{Tr}BA$ and $\operatorname{Tr}(A+B) = \operatorname{Tr}A + \operatorname{Tr}B$. Since $R_i\succcurlyeq 0$, $i\in\mathbb{N}_{[1,N^2]}$, we have also that
	\begin{equation*}
		Q\coloneqq C\sum_{\mathclap{i\in\mathbb{N}_{[1,N^2]}}}R_iC^\top \succcurlyeq 0,
	\end{equation*}
    It follows, due to~\eqref{eq:negsemdeftr} and $Q=Q^\top$, that $Q=0$. Using also that $R_{N^2+i}\succcurlyeq 0$, $i\in\mathcal{N}$, and $R_i\neq 0$ for some $i\in\mathbb{N}_{[1,N^2]}$ (otherwise $R=0$ by~\eqref{eq:lastRs}), we conclude that exactly one of the statements in Theorem~\ref{thm:toa-pqdet} holds. 
\end{proof}

\begin{proof}{Proof of Lemma~\ref{lem:interm-toa}}
    Consider the system~\eqref{eq:system}. Let $\mathcal{S}^{n}_N\subseteq\mathbb{S}^{Nn}$ denote, again, a space of block-diagonal Hermitian matrices, composed of $N$ diagonal blocks of size $n$-by-$n$, with inner product~\eqref{eq:inner-prod} and let $S\coloneqq \operatorname{diag}\{S_1,S_2,\hdots,S_{N}\}$ with symmetric matrices $S_i\in\mathbb{S}^{n_x}$, $i\in\mathcal{N}$. Then,~\eqref{eq:pqstab-cond} can be expressed as $\mathcal{A}(S)+A_0\succ 0$ with $\mathcal{A}\colon\mathcal{S}^{n_x}_{N}\rightarrow\mathcal{S}^{n_x}_{(N+1)N}$ (i.e., $\mathcal{V}\leftarrow\mathcal{S}^{n_x}_{N}$ and $\mathcal{S}\leftarrow\mathcal{S}^{n_x}_{(N+1)N}$ in Theorem~\ref{thm:thm1balakrishnan}) given by
    \begin{equation*}
        \mathcal{A}(S)\coloneqq \operatorname{diag}\{\operatorname{diag}\{S_j - A_iS_iA_i^\top\}_{i,j\in\mathcal{N}},S\}.
    \end{equation*}
    and $A_0\coloneqq \operatorname{diag}\{\operatorname{diag}\{BB^\top\}_{i,j\in\mathcal{N}},0_{Nn_x,Nn_x}\}$. Note that we use the convention $\{A_{ij}\}_{i,j\in\mathcal{N}}=\{A_{11},A_{21},A_{12},A_{22}\}$. Using the fact that $\operatorname{Tr}ABC=\operatorname{Tr}BCA$ and that $\mathcal{A}^{\operatorname{adj}}\colon\mathcal{S}^{n_x}_{(N+1)N}\rightarrow\mathcal{S}^{n_x}_{N}$ satisfies, for all $S_i\in\mathbb{S}^{n_x}$, $i\in\mathcal{N}$, and $R_i\in\mathbb{S}^{n_x}$, $i\in\mathbb{N}_{[1,(N+1)N]}$,
    \begin{align*}
        &\left\langle\mathcal{A}(S),R\right\rangle_{\mathcal{S}^{n_x}_{(N+1)N}}=\sum_{i\in\mathcal{N}}\operatorname{Tr}S_iR_{N^2+i} + \nonumber\\
        &\,\sum_{i\in\mathcal{N}}\sum_{j\in\mathcal{N}}\operatorname{Tr}(S_j-A_iS_iA_i^\top)R_{(j-1)N+i}=\nonumber\\
        &\sum_{i\in\mathcal{N}}\operatorname{Tr}S_i(R_{N^2+i}+\sum_{j\in\mathcal{N}}(R_{(i-1)N+j} - A_i^\top R_{(j-1)N+i}A_i))\nonumber\\
        &= \left\langle S,\mathcal{A}^{\operatorname{adj}}(R)\right\rangle_{S^{n_x}_{N}}
    \end{align*}
    with $R=\operatorname{diag}\{R_i\}_{i\in\mathbb{N}_{[1,(N+1)N]}}\in\mathcal{S}^{n_x}_{(N+1)N}$, we find
    \begin{align*}
        &\mathcal{A}^{\operatorname{adj}}(R) = \\
        &\operatorname{diag}\{R_{N^2+i}+\sum_{j\in\mathcal{N}}(R_{(i-1)N+j}-A_i^\top R_{(j-1)N+i}A_i)\}_{i\in\mathcal{N}}.\nonumber
    \end{align*}
    By Theorem~\ref{thm:thm1balakrishnan}, there exist $S_i\in\mathbb{S}^{n_x}$, $i\in\mathcal{N}$, satisfying the conditions in~\ref{item:interm-toa-1} if and only if there does not exist $R=\operatorname{diag}\{R_i\}_{i\in\mathbb{N}_{[1,(N+1)N]}}\in\mathcal{S}^{n_x}_{(N+1)N}$ such that $R\succcurlyeq 0$, $R\neq 0$ and $\mathcal{A}^{\operatorname{adj}}=0$ and $\left\langle A_0,R\right\rangle_{\mathcal{S}^{n_x}_{(N+1)N}}\leqslant 0$, where $\mathcal{A}^{\operatorname{adj}}(R)=0$ evaluates to, for all $i\in\mathcal{N}$,
    \begin{equation}
        R_{N^2+i} = \sum_{j\in\mathcal{N}}A_i^\top R_{(j-1)N+i}A_i - R_{(i-1)N+j},\label{eq:lastRs2}
    \end{equation}
    and we obtain
    \begin{equation}
        \left\langle A_0,R\right\rangle_{\mathcal{S}^{n_x}_{(N+1)N}}=\operatorname{Tr}B^\top\sum_{\mathclap{i\in\mathbb{N}_{[1,N^2]}}}R_iB\leqslant 0.
        \label{eq:Qpos}
    \end{equation}
    Since $R_i\succcurlyeq 0$, $i\in\mathbb{N}_{[1,N^2]}$, we have
    \begin{equation*}
        Q\coloneqq B^\top \sum_{\mathclap{i\in\mathbb{N}_{[1,N^2]}}}R_iB\succcurlyeq 0,
    \end{equation*}
    which, together with~\eqref{eq:Qpos}, implies that $Q=0$. This, together with the fact that $R_{N^2+i}\succcurlyeq 0$, $i\in\mathcal{N}$, and that $R_i\neq 0$ for some $i\in\mathbb{N}_{[1,N^2]}$ (since otherwise $R=0$ by~\eqref{eq:lastRs2}), leads us to conclude that exactly one of the statements~\ref{item:interm-toa-1} and~\ref{item:interm-toa-2} holds.  
\end{proof}

\begin{proof}{Proof of Theorem~\ref{thm:toa-pqstab}}
    Consider the system~\eqref{eq:system}, being strictly polytopic. By Theorem~\ref{thm:pQstab}, if the system is poly-Q stabilizable, then there exists $S_i\in\mathbb{S}^{n_x}_{\succ 0}$, $i\in\mathcal{N}$, satisfying~\eqref{eq:pqstab-cond} for all $i,j\in\mathcal{N}$. Hence,~\ref{item:toa-pqstab-1} implies~\ref{item:interm-toa-1}, such that, if~\ref{item:toa-pqstab-1} holds, then~\ref{item:toa-pqstab-2} does not hold. Thus, at most one of the statements in Theorem~\ref{thm:pQstab} holds.  
\end{proof}
                                       
\end{document}